\documentclass[12pt]{article}
\usepackage{latexsym, amssymb, amsmath, amscd, amsfonts, epsfig, graphicx, colordvi,verbatim,ifpdf,extarrows}
\usepackage{amsfonts, amsmath, amssymb}
\usepackage{amssymb,amsfonts,amsmath,latexsym,epsfig,cite, psfrag,eepic,color}
\usepackage{amscd,graphics}
\usepackage{latexsym, amssymb,  amsmath,amscd, amsfonts, epsfig, graphicx, colordvi,amsthm}

\usepackage{graphicx}
\usepackage{epstopdf}
\usepackage{color}
\usepackage{ifpdf}
\usepackage{fancybox}
\usepackage[font=small,labelfont=bf,labelsep=none]{caption}
\usepackage{float}

\allowdisplaybreaks

\newtheorem{thm}{Theorem}[section]

\newtheorem{defi}[thm]{Definition}
\newtheorem{lem}[thm]{Lemma}
\newtheorem{core}[thm]{Corollary}

\def\pf{\noindent{\it Proof.} }
\def\qed{\nopagebreak\hfill{\rule{4pt}{7pt}}
\medbreak}

\numberwithin{equation}{section}

\def\qed{\nopagebreak\hfill{\rule{4pt}{7pt}}
\medbreak}

\newlength{\boxedparwidth}
  {\begin{center} \begin{tabular}{|@{\hspace{.315in}}c@{\hspace{.15in}}|}
                  \hline \\ \begin{minipage}[t]{\boxedparwidth}
                  \setlength{\parindent}{.25in}}%
  {\end{minipage} \\ \\ \hline \end{tabular} \end{center}}

\begin{document}

\begin{center}
{\Large \bf Parity considerations in the number of parts}
\end{center}

\begin{center}
{Thomas Y. He}$^{1}$, {H.X. Huang}$^{2}$,
  {Y.X. Xie}$^{3}$ and {T.T. Zou}$^{4}$ \vskip 2mm

$^{1,2,3,4}$ School of Mathematical Sciences, Sichuan Normal University, Chengdu 610066, P.R. China

   \vskip 2mm

  $^1$heyao@sicnu.edu.cn, $^2$huanghaoxuan@stu.sicnu.edu.cn,  $^3$xieyx@stu.sicnu.edu.cn, $^4$zoutingting@stu.sicnu.edu.cn
\end{center}

\vskip 6mm   {\noindent \bf Abstract.} Recently, Chen, He, Hu and Xie considered the parity of the number of non-overlined (resp. overlined)  parts of size greater than or equal to the size of the smallest overlined (resp. non-overlined) part  in an overpartition. In this article,  we investigate the parity of the number of non-overlined (resp. overlined)  parts of size less than or equal to the size of the largest overlined (resp. non-overlined) part  in an overpartition. We also study the parity of the number of even parts greater than the smallest odd part  in a partition.



\section{Introduction}

A partition $\pi$ of a positive integer $n$ is a finite non-increasing sequence of positive integers $\pi=(\pi_1,\pi_2,\ldots,\pi_m)$ such that $\pi_1+\pi_2+\cdots+\pi_m=n$. The empty sequence forms the only partition of zero. The $\pi_i$ are called the parts of $\pi$. Let $\ell(\pi)$ be the number of parts of $\pi$.

In \cite{Andrews-2018,Andrews-2019}, Andrews considered partitions in which parts of a given parity are all smaller than those of the other parity, and if the smaller parity is odd then odd parts must appear. One of the objectives of this article is to study ${p}^{uu}(n)$, which is the number of partitions of $n$ such that odd parts must appear and each even part is smaller than each odd part.
Moreover, we introduce the following three partition functions.

\begin{itemize}
\item Let ${p}^{ud}(n)$ be the number of partitions counted by ${p}^{uu}(n)$ with distinct odd parts.

\item Let ${p}^{du}(n)$ be the number of partitions counted by ${p}^{uu}(n)$ with distinct even parts.

\item Let ${p}^{dd}(n)$ be the number of partitions counted by ${p}^{uu}(n)$ with distinct parts.
\end{itemize}

We will give the generating functions of ${p}^{uu}(n)$, ${p}^{ud}(n)$,  ${p}^{du}(n)$ and ${p}^{dd}(n)$.

\begin{thm}\label{gen-e<o}
\begin{equation}\label{gen-e<o-uu}
\sum_{n\geq 1}{p}^{uu}(n)q^n=\frac{q}{(1-q)(q^2;q^2)_\infty},
\end{equation}
\begin{equation}\label{gen-e<o-ud}
\sum_{n\geq 1}{p}^{ud}(n)q^n=\frac{1}{(q^2;q^2)_\infty}\sum_{k=1}^\infty q^{k^2},
\end{equation}
\begin{equation}\label{gen-e<o-du}
\sum_{n\geq 1}{p}^{du}(n)q^n=\sum_{m\geq 1}\frac{q^{m^2+m}}{(q^2;q^2)_m}\sum_{j=1}^{m}q^{-j^2},
\end{equation}
\begin{equation}\label{gen-e<o-dd}
\sum_{n\geq 1}{p}^{dd}(n)q^n=\frac{1}{1-q}\left((-q;q^2)_\infty-(-q^2;q^2)_\infty\right).
\end{equation}
\end{thm}

Inspired by the partitions with parts separated by parity, Chen, He, Hu and Xie \cite{Chen-He-Hu-Xie-2024} considered the overpartitions with separated overlined parts and non-overlined parts.
An overpartition, introduced by Corteel and Lovejoy \cite{Corteel-Lovejoy-2004},  is a partition such that the first occurrence of a number can be overlined. For example, there are eight overparitions of $3$.
\[(3),(\overline{3}),(2,1),(\overline{2},1),(2,\overline{1}),(\overline{2},\overline{1}),(1,1,1),(\overline{1},1,1).\]

For a part $\pi_i$ of $\pi$, we say that $\pi_i$ is of size $t$ if $\pi_i=t$ or $\overline{t}$. For easier expression, Chen, He, Hu and Xie \cite{Chen-He-Hu-Xie-2024} introduced the following notations.
\begin{itemize}
\item Let $LN(\pi)$ (resp. $SN(\pi)$) be the size of the largest (resp. smallest) non-overlined part of $\pi$ if there exist non-overlined parts in $\pi$, and $LN(\pi)=0$ (resp. $SN(\pi)=0$) otherwise.

\item Let $LO(\pi)$ (resp. $SO(\pi)$) be the size of the largest (resp. smallest) overlined part of $\pi$ if there exist overlined parts in $\pi$, and $LO(\pi)=0$ (resp. $SO(\pi)=0$) otherwise.

\item Let $\ell_{N\geq O}(\pi)$ (resp. $\ell_{N>O}(\pi)$) be the number of non-overlined parts  with size greater than or equal to (resp. greater than) $SO(\pi)$ in $\pi$.

\item Let $\ell_{O\geq N}(\pi)$ (resp. $\ell_{O>N}(\pi)$) be the number of overlined parts of size greater than or equal to (resp. greater than) $SN(\pi)$ in $\pi$.

\end{itemize}

Motivated by Kim, Kim and Lovejoy \cite{Kim-Kim-Lovejoy-2021} and Lin and Lin \cite{Lin-Lin-2024},  Chen, He, Hu and Xie \cite{Chen-He-Hu-Xie-2024} considered the following partition functions.
 \begin{itemize}
\item[(1)] Let $A_{N\geq O}(n)$ (resp. $B_{N\geq O}(n)$) be the number of overpartitions $\pi$ of $n$ with $\ell_{N\geq O}(\pi)$ being even (resp. odd).

\item[(2)] Let $A_{N>O}(n)$ (resp. $B_{N>O}(n)$) be the number of overpartitions $\pi$ of $n$ with $\ell_{N>O}(\pi)$ being even (resp. odd).

\item[(3)] Let $A_{O\geq N}(n)$ (resp. $B_{O\geq N}(n)$) be the number of overpartitions $\pi$ of $n$ such that
 $SN(\pi)\geq 1$ and $\ell_{O\geq N}(\pi)$ is even (resp. odd).

\item[(4)] Let $A_{O>N}(n)$ (resp. $B_{O>N}(n)$) be the number of overpartitions $\pi$ of $n$ such that $SN(\pi)\geq 1$ and $\ell_{O>N}(\pi)$ is even (resp. odd).
\end{itemize}

 Chen, He, Hu and Xie \cite{Chen-He-Hu-Xie-2024} obtained the following identities.
\begin{thm}
For $n\geq1$,
\begin{equation}\label{lem-excess-old-1}
A_{N\geq O}(n)-B_{N\geq O}(n)=2\left(p^{e}_o(n)-p^{o}_e(n)\right),
\end{equation}
where $p^{e}_o(n)$ {\rm(}resp. $p^{o}_e(n)${\rm)} is the number of partitions of $n$ such that the largest part is even {\rm(}resp. odd{\rm)} and the smallest part is odd {\rm(}resp. even{\rm);}
\begin{equation}\label{EE-OE-eqn-1}
A_{N> O}(n)-B_{N> O}(n)=2p_e(n),
\end{equation}
where $p_e(n)$ is the number of partitions of $n$ with an even number of parts{\rm;}
\begin{equation}\label{EG-OG-eqn-D}
A_{O\geq N}(n)-B_{O\geq N}(n)=D(n),
\end{equation}
where $D(n)$ is the number of distinct partitions of $n${\rm;}
\[A_{O> N}(n)-B_{O> N}(n)={H}'_{ON}(n),\]
where ${H}'_{ON}(n)$ is the number of overpartitions $\pi$ of $n$ with $LN(\pi)=SN(\pi)\geq 1$ and $SN(\pi)\geq LO(\pi)$.
\end{thm}

In this article, we first study the remaining four cases.
\begin{itemize}
\item[(1)] Let $\overline{A}_{N\geq O}(n)$ (resp. $\overline{B}_{N\geq O}(n)$) be the number of overpartitions $\pi$ of $n$ such that $SO(\pi)\geq 1$ and $\ell_{N\geq O}(\pi)$ is even (resp. odd).

\item[(2)] Let $\overline{A}_{N>O}(n)$  (resp. $\overline{B}_{N>O}(n)$) be the number of overpartitions $\pi$ of $n$ such that $SO(\pi)\geq 1$ and $\ell_{N>O}(\pi)$ is even (resp. odd).

\item[(3)] Let $\overline{A}_{O\geq N}(n)$ (resp. $\overline{B}_{O\geq N}(n)$) be the number of overpartitions $\pi$ of $n$ with $\ell_{O\geq N}(\pi)$ being even (resp. odd).

\item[(4)] Let $\overline{A}_{O>N}(n)$ (resp. $\overline{B}_{O>N}(n)$) be the number of overpartitions $\pi$ of $n$ with $\ell_{O>N}(\pi)$ being even (resp. odd).
\end{itemize}

 With a similar argument in \cite{Chen-He-Hu-Xie-2024}, we get the following identities.
\begin{thm}\label{remin-thm-over}
 For $n\geq1$,
\begin{equation}\label{lem-remain-excess-1}
\overline{A}_{N\geq O}(n)-\overline{B}_{N\geq O}(n)=p^{o}(n)-p^{e}(n),
\end{equation}
where $p^{o}(n)$ {\rm(}resp. $p^{e}(n)${\rm)} is the number of partitions of $n$ such that the largest part appears an odd {\rm(}resp. even{\rm)} number of times{\rm;}
\begin{equation}\label{EE-OE-remain-eqn-1}
\overline{A}_{N>O}(n)-\overline{B}_{N>O}(n)=p(n),
\end{equation}
where $p(n)$ is the number of partitions of $n${\rm;}
\begin{equation}\label{remian-eqn-3}
\overline{A}_{O\geq N}(n)-\overline{B}_{O\geq N}(n)=2D_{e}(n),
\end{equation}
where $D_{e}(n)$ is the number of distinct partitions of $n$ with an even number of parts{\rm;}
\begin{equation}\label{remian-eqn-4}
\overline{A}_{O>N}(n)-\overline{B}_{O>N}(n)=2{H}^{o}_{ON}(n),
\end{equation}
where ${H}^{o}_{ON}(n)$ is the number of overpartitions $\pi$ of $n$ such that $LN(\pi)=SN(\pi)\geq 1$, $SN(\pi)\geq LO(\pi)$ and there is  an odd number of overlined parts.
\end{thm}

As a corollary, we can get the following result.
\begin{core}\label{remin-core-over}
For $n\geq 1$,
\begin{equation}\label{core-remain-excess-4}
\overline{A}_{N\geq O}(n)-\overline{B}_{N\geq O}(n)\geq 0\text{ with strict inequality if }n\neq 2,
\end{equation}
\begin{equation}\label{core-remain-excess-1}
\overline{A}_{N>O}(n)-\overline{B}_{N>O}(n)>0,
\end{equation}
\begin{equation}\label{core-remain-excess-2}
\overline{A}_{O\geq N}(n)-\overline{B}_{O\geq N}(n)\geq 0\text{ with strict inequality if }n\geq 3,
\end{equation}
\begin{equation}\label{core-remain-excess-3}
\overline{A}_{O>N}(n)-\overline{B}_{O>N}(n)\geq 0\text{ with strict inequality if }n\geq 2.
\end{equation}
\end{core}

Then, we will investigate the number of non-overlined (resp. overlined)  parts of size less than or equal to the size of the largest overlined (resp. non-overlined) part in an overpartition.
For easier expression, we introduce the following notations.
\begin{itemize}
\item Set $\widetilde{LN}(\pi)={LN}(\pi)$ and $\widetilde{SN}(\pi)={SN}(\pi)$ if there exist non-overlined parts in $\pi$, and set $\widetilde{LN}(\pi)=+\infty$ and $\widetilde{SN}(\pi)=+\infty$ otherwise.

\item Set $\widetilde{LO}(\pi)={LO}(\pi)$ and $\widetilde{SO}(\pi)={SO}(\pi)$ if there exist overlined parts in $\pi$, and set $\widetilde{LO}(\pi)=+\infty$ and $\widetilde{SO}(\pi)=+\infty$ otherwise.

\item Let $\ell_{N\leq O}(\pi)$ (resp. $\ell_{N<O}(\pi)$) be the number of non-overlined parts  of size less than or equal to (resp. less than) $\widetilde{LO}(\pi)$ in $\pi$.

\item Let $\ell_{O\leq N}(\pi)$ (resp. $\ell_{O<N}(\pi)$) be the number of overlined parts of size less than or equal to (resp. less than) $\widetilde{LN}(\pi)$ in $\pi$.

\end{itemize}

For example, the notations above  for the overpartitions of $3$ are given in the following table.
{ \begin{center}
\begin{tabular}{|c|c|c|c|c|c|c|c|c|}
  \hline
  $\pi$&$\widetilde{LN}(\pi)$&$\widetilde{SN}(\pi)$&$\widetilde{LO}(\pi)$&$\widetilde{SO}(\pi)$&$\ell_{N\leq O}(\pi)$&$\ell_{N<O}(\pi)$&$\ell_{O\leq N}(\pi)$&$\ell_{O<N}(\pi)$\\
  \hline
  $(3)$&$3$&$3$&$+\infty$&$+\infty$&$1$&$1$&$0$&$0$\\
  \hline
  $(\bar{3})$&$+\infty$&$+\infty$&$3$&$3$&$0$&$0$&$1$&$1$\\
  \hline
  $(2,1)$&$2$&$1$&$+\infty$&$+\infty$&$2$&$2$&$0$&$0$\\
  \hline
  $(\bar{2},1)$&$1$&$1$&$2$&$2$&$1$&$1$&$0$&$0$\\
  \hline
  $(2,\bar{1})$&$2$&$2$&$1$&$1$&$0$&$0$&$1$&$1$\\
  \hline
  $(\bar{2},\bar{1})$&$+\infty$&$+\infty$&$2$&$1$&$0$&$0$&$2$&$2$\\
  \hline
  $(1,1,1)$&$1$&$1$&$+\infty$&$+\infty$&$3$&$3$&$0$&$0$\\
  \hline
  $(\bar{1},1,1)$&$1$&$1$&$1$&$1$&$2$&$0$&$1$&$0$\\
  \hline
\end{tabular}
\end{center}}

We consider the following partition functions.

\begin{itemize}
\item[(1)] Let $A_{O\leq N}(n)$ (resp. $B_{O\leq N}(n)$) be the number of overpartitions $\pi$ of $n$ with $\ell_{O\leq N}(\pi)$ being even (resp. odd).

\item[(2)] Let $\widetilde{A}_{O\leq N}(n)$ (resp. $\widetilde{B}_{O\leq N}(n)$) be the number of overpartitions $\pi$ of $n$ such that non-overlined parts appear  and $\ell_{O\leq N}(\pi)$ is even (resp. odd).

\item[(3)] Let $A_{O<N}(n)$ (resp. $B_{O<N}(n)$) be the number of overpartitions $\pi$ of $n$ with $\ell_{O<N}(\pi)$ being even (resp. odd).

\item[(4)] Let $\widetilde{A}_{O<N}(n)$ (resp. $\widetilde{B}_{O<N}(n)$) be the number of overpartitions $\pi$ of $n$ such that non-overlined parts appear  and $\ell_{O<N}(\pi)$ is even (resp. odd).

\item[(5)] Let $A_{N\leq O}(n)$ (resp. $B_{N\leq O}(n)$) be the number of overpartitions $\pi$ of $n$ with $\ell_{N\leq O}(\pi)$ being even (resp. odd).

\item[(6)] Let $\widetilde{A}_{N\leq O}(n)$ (resp. $\widetilde{B}_{N\leq O}(n)$) be the number of overpartitions $\pi$ of $n$ such that overlined parts appear  and $\ell_{N\leq O}(\pi)$ is even (resp. odd).

\item[(7)] Let $A_{N<O}(n)$ (resp. $B_{N<O}(n)$) be the number of overpartitions $\pi$ of $n$ with $\ell_{N<O}(\pi)$ being even (resp. odd).

\item[(8)] Let $\widetilde{A}_{N<O}(n)$ (resp. $\widetilde{B}_{N<O}(n)$) be the number of overpartitions $\pi$ of $n$ such that overlined parts appear  and $\ell_{N<O}(\pi)$ is even (resp. odd).

\end{itemize}

For example, the partition functions above  for $n=3$ are given in the following tables.
\begin{center}
\begin{tabular}{|c|c|c|c|c|c|c|c|}
  \hline
  $A_{O\leq N}(3)$&$B_{O\leq N}(3)$&$\widetilde{A}_{O\leq N}(3)$&$\widetilde{B}_{O\leq N}(3)$&$A_{O<N}(3)$&$B_{O<N}(3)$&$\widetilde{A}_{O<N}(3)$&$\widetilde{B}_{O<N}(3)$\\
  \hline
  $5$&$3$&$4$&$2$&$6$&$2$&$5$&$1$\\
  \hline
\end{tabular}
\end{center}

\begin{center}
\begin{tabular}{|c|c|c|c|c|c|c|c|c|}
  \hline
  $A_{N\leq O}(3)$&$B_{N\leq O}(3)$&$\widetilde{A}_{N\leq O}(3)$&$\widetilde{B}_{N\leq O}(3)$&$A_{N<O}(3)$&$B_{N<O}(3)$&$\widetilde{A}_{N<O}(3)$&$\widetilde{B}_{N<O}(3)$\\
  \hline
  $5$&$3$&$4$&$1$&$5$&$3$&$4$&$1$\\
  \hline
\end{tabular}
\end{center}

We get the following identities.
\begin{thm}\label{new-A-B-thm}
For $n\geq1$,
\begin{equation}\label{E-O-EQN(1)}
A_{O\leq N}(n)-B_{O\leq N}(n)=2D_e(n),
\end{equation}
where $D_{e}(n)$ is the number of distinct partitions of $n$ with an even number of parts{\rm;}
\begin{equation}\label{E-O-EQN(2)}
\widetilde{A}_{O\leq N}(n)-\widetilde{B}_{O\leq N}(n)=D(n),
\end{equation}
where $D(n)$ is the number of distinct partitions of $n${\rm;}
\begin{equation}\label{E-O-EQN(3)}
A_{O<N}(n)-B_{O<N}(n)=2H_{N<O}^o(n),
\end{equation}
where $H_{N<O}^o(n)$ is the number of overpartitions $\pi$ of $n$ such that non-overlined parts appear, $\widetilde{SO}(\pi)\geq \widetilde{LN}(\pi)=\widetilde{SN}(\pi)$ and there is an odd number of overlined parts{\rm;}
\begin{equation}\label{E-O-EQN(4)}
\widetilde{A}_{O<N}(n)-\widetilde{B}_{O<N}(n)=H_{N<O}(n),
\end{equation}
where $H_{N<O}(n)$ is the number of overpartitions $\pi$ of $n$ such that non-overlined parts appear and $\widetilde{SO}(\pi)\geq \widetilde{LN}(\pi)=\widetilde{SN}(\pi)${\rm;}
\begin{equation}\label{E-O-EQN(5)}
A_{N\leq O}(n)-B_{N\leq O}(n)=2\left(p''_{o}(n)-p''_{e}(n)\right),
\end{equation}
where $p''_{o}(n)$ {\rm(}resp. $p''_{e}(n)${\rm)} is the number of partitions $\pi$ of $n$ such that $\widetilde{SN}(\pi)$ appears an odd {\rm(}resp. even{\rm)} number of times and there is an odd number of parts greater than $\widetilde{SN}(\pi)${\rm;}
\begin{equation}\label{E-O-EQN(6)}
\widetilde{A}_{N\leq O}(n)-\widetilde{B}_{N\leq O}(n)=p'_o(n)-p'_e(n),
\end{equation}
where $p'_o(n)$ {\rm(}resp. $p'_e(n)${\rm)} is the number of partitions $\pi$ of $n$ such that $\widetilde{SN}(\pi)$ appears an odd {\rm(}resp. even{\rm)} number of times{\rm;}
\begin{equation}\label{E-O-EQN(7)}
A_{N<O}(n)-B_{N<O}(n)=2p_e(n),
\end{equation}
where $p_e(n)$ is the number of partitions of $n$ with an even number of parts{\rm;}
\begin{equation}\label{E-O-EQN(8)}
\widetilde{A}_{N<O}(n)-\widetilde{B}_{N<O}(n)=p(n),
\end{equation}
where $p(n)$ is the number of partitions of $n${\rm.}
\end{thm}

By \eqref{EE-OE-eqn-1}, \eqref{EG-OG-eqn-D}, \eqref{EE-OE-remain-eqn-1}, \eqref{remian-eqn-3}, \eqref{E-O-EQN(1)}, \eqref{E-O-EQN(2)}, \eqref{E-O-EQN(7)} and \eqref{E-O-EQN(8)}, we can get the following two corollaries.
\begin{core}\label{coro-equal-not}
For $n\geq 1$,
\[\overline{A}_{O\geq N}(n)={A}_{O\leq N}(n),\]
\[\overline{B}_{O\geq N}(n)={B}_{O\leq N}(n),\]
\[{A}_{O\geq N}(n)=\widetilde{A}_{O\leq N}(n),\]
\[{B}_{O\geq N}(n)=\widetilde{B}_{O\leq N}(n).\]
\end{core}

\begin{core}\label{coro-greater-not}
For $n\geq 1$,
\[{A}_{N>O}(n)={A}_{N<O}(n),\]
\[{B}_{N>O}(n)={B}_{N<O}(n),\]
\[\overline{A}_{N>O}(n)=\widetilde{A}_{N<O}(n),\]
\[\overline{B}_{N>O}(n)=\widetilde{B}_{N<O}(n).\]
\end{core}

We will give a combinatorial proof of Corollary \ref{coro-equal-not}. It would be interesting to give a combinatorial proof of Corollary \ref{coro-greater-not}.

Let $\pi$ be a partition such that there exist odd parts in $\pi$. We use $\ell_{e>o}(\pi)$ to denote the number of even parts greater than the smallest odd part in $\pi$. It is clear that if $\pi$ is a partition enumerated by ${p}^{uu}(n)$ then $\ell_{e>o}(\pi)=0$. Then, we introduce the following partition functions.
 \begin{itemize}
 \item[(1)] Let ${p}^{uu}_e(n)$ (resp. ${p}^{uu}_o(n)$) be the number of partitions $\pi$ of $n$ such that odd parts appear in $\pi$ and $\ell_{e>o}(\pi)$ is even (resp. odd).

  \item[(2)] Let ${p}^{ud}_e(n)$ (resp. ${p}^{ud}_o(n)$) be the number of partitions $\pi$ of $n$ such that odd parts appear in $\pi$, odd parts are distinct and $\ell_{e>o}(\pi)$ is even (resp. odd).

       \item[(3)] Let ${p}^{du}_e(n)$ (resp. ${p}^{du}_o(n)$) be the number of partitions $\pi$ of $n$ such that odd parts appear in $\pi$, even parts are distinct and $\ell_{e>o}(\pi)$ is even (resp. odd).

        \item[(4)] Let ${p}^{dd}_e(n)$ (resp. ${p}^{dd}_o(n)$) be the number of partitions $\pi$ of $n$ such that odd parts appear in $\pi$, the parts are distinct and $\ell_{e>o}(\pi)$ is even (resp. odd).
 \end{itemize}

 We have the following results.
 \begin{thm}\label{thm-eo-excess}
 For $n\geq 1$,
 \begin{itemize}
 \item[\rm{(1)}] ${p}^{uu}_e(n)-{p}^{uu}_o(n)$ equals the number of pairs of partitions $(\alpha,\beta)$ such that the sum of the parts in $\alpha$ and $\beta$ is $n$, $\alpha$ consists of distinct even parts, the even parts of $\beta$ are divisible by $4$, and odd parts appear in $\beta$ which are equal to $2\ell(\alpha)+1$;

  \item[\rm{(2)}] ${p}^{ud}_e(n)=\overline{p}^{ud}_e(n)$ and ${p}^{ud}_o(n)=\overline{p}^{ud}_o(n)$, where $\overline{p}^{ud}_e(n)$ (resp. $\overline{p}^{ud}_o(n)$) is the number of pairs of partitions $(\alpha,\beta)$ such that the sum of the parts in $\alpha$ and $\beta$ is $n$, the parts of $\alpha$ are even, $\beta$ is nonempty, $\beta$ consists of distinct odd parts, and $\frac{\beta_1+1}{2}+\ell(\beta)$ is even (resp. odd);

  \item[\rm{(3)}]   ${p}^{du}_e(n)-{p}^{du}_o(n)$ equals the number of partitions $\pi=(\pi_1,\pi_2,\ldots,\pi_m)$ of $n$ such that $\pi_1,\pi_2,\ldots,\pi_{m-1}$ are distinct even parts which are greater than $2\pi_m$;

  \item[\rm{(4)}] ${p}^{dd}_e(n)-{p}^{dd}_o(n)$ equals the number of partitions of $n$  such that odd parts must appear,  each positive odd integer not exceeding the largest odd part  appears once, even parts are distinct  and each even part is greater than one plus the largest odd part.

 \end{itemize}
 \end{thm}

This article is organized as follows.  We collect some necessary notations and results in Section 2.
 We  prove Theorem \ref{gen-e<o} in Section 3 and Theorem \ref{thm-eo-excess} in Section 4. The proofs of Theorem \ref{remin-thm-over} and Corollary \ref{remin-core-over} are presented in Section 5. We give an analytic proof and a combinatorial proof of Theorem \ref{new-A-B-thm} in Section 6 and Section 7 respectively.  A combinatorial proof of Corollary \ref{coro-equal-not} is given in Section 8.

\section{Preliminaries}

In this section, we collect some notations and results needed in this article.
In the rest of this article, we use the following notations.
\begin{itemize}
\item $\ell_o(\pi)$: the number of overlined parts in $\pi$.

\item ${\mathcal{P}}(n)$: the set of partitions of $n$.

\item ${\mathcal{D}}(n)$: the set of distinct partitions of $n$.

\item $\overline{\mathcal{P}}(n)$: the set of overpartitions of $n$.

\item $\overline{\mathcal{D}}(n)$: the set of overpartitions of $n$ with no non-overlined parts.

\item $\overline{\mathcal{PN}}(n)$: the set of overpartitions of $n$ such that non-overlined parts appear.

\item $\overline{\mathcal{PO}}(n)$: the set of overpartitions of $n$ such that overlined parts appear.
\end{itemize}

 We assume that $|q|<1$ and use the standard notation \cite{Andrews-1976}:
\[(a;q)_\infty=\prod_{i=0}^{\infty}(1-aq^i)\quad\text{and}\quad(a;q)_n=\frac{(a;q)_\infty}{(aq^n;q)_\infty}.\]

{\noindent \bf The $q$-binomial theorem \cite[Theorem 2.1]{Andrews-1976}:}
\begin{equation}\label{q-bin-eqn}
\sum_{n\geq 0}\frac{(a;q)_{n}}{(q;q)_{n}}t^{n}=\frac{{(at;q)}_{\infty}}{(t;q)_{\infty}}.
\end{equation}

{\noindent \bf A formula due to Euler \cite{Euler-1748} (see also \cite[(2.2.6)]{Andrews-1976}):}
\begin{equation}\label{Euler-2}
\sum_{n\geq0}\frac{t^nq^{{n}\choose 2}}{(q;q)_n}=(-t;q)_\infty.
\end{equation}

{\noindent \bf Jacobi's triple product identity \cite{Jacobi-1829} (see also \cite[Theorem 2.8]{Andrews-1976}):}
\begin{equation}\label{jacobi-triple}
\sum_{n=-\infty}^{\infty}z^nq^{n^2}=(-zq;q^2)_\infty(-q/z;q^2)_\infty(q^2;q^2)_\infty.
\end{equation}

{\noindent \bf Heine's transformation formula \cite{Heine-1847} (see also \cite[Corollary 2.3]{Andrews-1976}):}
\begin{equation}\label{Heine}
\sum_{n\geq 0}\frac{(a;q)_n(b;q)_n}{(q;q)_n(c;q)_n}t^n=\frac{(b;q)_\infty(at;q)_\infty}{(c;q)_\infty(t;q)_\infty}\sum_{n\geq 0}\frac{(c/b;q)_n(t;q)_n}{(q;q)_n(at;q)_n}b^n.
\end{equation}

By considering the number of parts, we have
\begin{equation}\label{gen-p}
\sum_{n\geq 0}q^n\sum_{\pi\in\mathcal{P}(n)}t^{\ell(\pi)}=\frac{1}{(tq;q)_\infty}=\sum_{n\geq 0}\frac{t^nq^n}{(q;q)_n}=1+\sum_{n\geq 1}\frac{tq^n}{(tq^{n};q)_\infty},
\end{equation}
and
\begin{equation}\label{gen-d}
\sum_{n\geq 0}q^n\sum_{\pi\in\mathcal{D}(n)}t^{\ell(\pi)}=(-tq;q)_\infty=1+\sum_{n\geq 1}tq^n(-tq;q)_{n-1}=1+\sum_{n\geq 1}tq^n(-tq^{n+1};q)_\infty.
\end{equation}

Setting $t=1$ in \eqref{gen-p} and \eqref{gen-d}, we get
\begin{equation}\label{gen-p+}
\sum_{n\geq0}p(n)q^{n}=\frac{1}{(q;q)_\infty}=\sum_{n\geq 0}\frac{q^n}{(q;q)_n}=1+\sum_{n\geq 1}\frac{q^n}{(q^{n};q)_\infty},
\end{equation}
and
\begin{equation}\label{gen-d+}
\sum_{n\geq0}D(n)q^{n}=(-q;q)_\infty=1+\sum_{n\geq 1}q^n(-q;q)_{n-1}=1+\sum_{n\geq 1}q^n(-q^{n+1};q)_\infty.
\end{equation}

Setting $t=-1$ in \eqref{gen-p} and \eqref{gen-d}, we get
\begin{equation}\label{gen-p-}
\sum_{n\geq0}\left(p_e(n)-p_o(n)\right)q^{n}=\frac{1}{(-q;q)_\infty}=\sum_{n\geq 0}\frac{(-1)^nq^n}{(q;q)_n}=1-\sum_{n\geq 1}\frac{q^n}{(-q^{n};q)_\infty},
\end{equation}
and
\begin{equation}\label{gen-d-}
\sum_{n\geq0}\left(D_e(n)-D_o(n)\right)q^{n}=(q;q)_\infty=1-\sum_{n\geq 1}q^n(q;q)_{n-1}=1-\sum_{n\geq 1}q^n(q^{n+1};q)_\infty,
\end{equation}
where $p_o(n)$ (resp. $D_o(n)$) is the number of partitions (resp. distinct partitions) of $n$ with an odd number of parts.

We conclude this section with an involution $\varphi$ on $\overline{\mathcal{P}}(n)$.
\begin{defi}\label{defi-involution-varphi}
For $n\geq 1$, let $\pi$ be an overpartition  in $\overline{\mathcal{P}}(n)$. The map $\varphi$
is defined as follows:
\begin{itemize}
\item[{\rm(1)}] if $\widetilde{SN}(\pi)<\widetilde{SO}(\pi)$, then change the smallest non-overlined part of $\pi$ to an overlined part{\rm;}

    \item[{\rm(2)}] if $\widetilde{SO}(\pi)\leq \widetilde{SN}(\pi)$, then change the smallest overlined part of $\pi$ to a non-overlined part{\rm.}
\end{itemize}
\end{defi}

\section{Proof of Theorem \ref{gen-e<o}}

The objective of this article is to give a proof of Theorem \ref{gen-e<o}. More precisely, we will show \eqref{gen-e<o} and \eqref{gen-e<o-ud} in Section 3.1 and  \eqref{gen-e<o-du} and \eqref{gen-e<o-dd} in Section 3.2.

\subsection{Proofs of \eqref{gen-e<o} and \eqref{gen-e<o-ud}}

With a similar argument in \cite[Section 2]{Andrews-2018}, we first give a proof of \eqref{gen-e<o}.

{\noindent \bf Proof of \eqref{gen-e<o}.} By considering the smallest odd part, we have
\begin{align*}
\sum_{n\geq 1}{p}^{uu}(n)q^n&=\sum_{n\geq 0}\frac{q^{2n+1}}{(q^{2n+1};q^2)_\infty(q^2;q^2)_n}\\
&=\frac{q}{(q;q^2)_\infty}\sum_{n\geq 0}\frac{(q;q^2)_nq^{2n}}{(q^2;q^2)_{n}}\\
&=\frac{q}{(q;q^2)_\infty}\frac{(q^3;q^2)_\infty}{(q^2;q^2)_\infty}\\
&=\frac{q}{(1-q)(q^2;q^2)_\infty},
\end{align*}
where the second equation follows from \eqref{q-bin-eqn} with $q\rightarrow q^2$, $a=q$ and $t=q^2$. This completes the proof.  \qed

Then, we give a proof of \eqref{gen-e<o-ud} with a similar argument in \cite[Section 3]{Andrews-2019}.

{\noindent \bf Proof of \eqref{gen-e<o-ud}.} In terms of the smallest odd part, we get
\begin{align*}
\sum_{n\geq 1}{p}^{ud}(n)q^n&=\sum_{n\geq 0}q^{2n+1}(-q^{2n+3};q^2)_\infty\frac{1}{(q^2;q^2)_n}\\
&=(-q;q^2)_\infty\sum_{n\geq 0}\frac{q^{2n+1}}{(-q;q^2)_{n+1}(q^2;q^2)_n}\\
&=(-q;q^2)_\infty\sum_{n\geq 0}\frac{q^{2n+1}}{(-q;-q)_{2n+1}}\\
&=(-q;q^2)_\infty\sum_{n\geq 0}\frac{q^n}{(-q;-q)_n}\frac{1-(-1)^n}{2}\\
&=\frac{(-q;q^2)_\infty}{2}\left(\sum_{n\geq 0}\frac{q^n}{(-q;-q)_n}-\sum_{n\geq 0}\frac{(-q)^n}{(-q;-q)_n}\right)\\
&=\frac{(-q;q^2)_\infty}{2}\left(\frac{1}{(q;-q)_\infty}-\frac{1}{(-q;-q)_\infty}\right)\\
&=\frac{1}{2(q^2;q^2)_\infty}\left((-q;q^2)_\infty^2(q^2;q^2)_\infty-1\right)\\
&=\frac{1}{2(q^2;q^2)_\infty}\left(\sum_{k=-\infty}^\infty q^{k^2}-1\right)\\
&=\frac{1}{(q^2;q^2)_\infty}\sum_{k=1}^\infty q^{k^2},
\end{align*}
where the last forth equation follows from \eqref{q-bin-eqn} and the last second equation follows from \eqref{jacobi-triple}. The proof  is complete.  \qed

\subsection{Proofs of \eqref{gen-e<o-du} and \eqref{gen-e<o-dd}}

Recently, Passary \cite[Section 3.2]{Passary-2019} and Chen, He, Tang and Wei  \cite[Section 3]{Chen-He-Tang-Wei-2024} studied the partitions with parts separated by parity in view of separable integer partition classes introduced by Andrews \cite{Andrews-2022}. We first give a proof of \eqref{gen-e<o-du} with a similar argument in \cite[Section 3.3]{Chen-He-Tang-Wei-2024}.

{\noindent \bf Proof of \eqref{gen-e<o-du}.} It is easy to see that
\begin{equation*}
\sum_{n\geq 1}{p}^{du}(n)q^n=\sum_{m\geq 1}\frac{1}{(q^2;q^2)_m}\sum_{j=1}^{m}q^{2+4+\cdots+(2m-2j)+j(2m-2j+1)}=\sum_{m\geq 1}\frac{1}{(q^2;q^2)_m}\sum_{j=1}^{m}q^{m^2+m-j^2}.
\end{equation*}
This completes the proof.  \qed

Then, we give a proof of \eqref{gen-e<o-dd} with a similar argument in \cite[Section 3.1]{Chen-He-Tang-Wei-2024}.

{\noindent \bf Proof of \eqref{gen-e<o-dd}.} It is clear that
\begin{align*}
\sum_{n\geq 1}{p}^{dd}(n)q^n&=\sum_{m\geq 1}\frac{1}{(q^2;q^2)_m}\sum_{j=0}^{m-1}q^{2+4+\cdots+2j+(2j+1)+\cdots+(2m-1)}\\
&=\sum_{m\geq 1}\frac{1}{(q^2;q^2)_m}\sum_{j=0}^{m-1}q^{m^2+j}\\
&=\sum_{m\geq 1}\frac{q^{m^2}}{(q^2;q^2)_m}\frac{1-q^m}{1-q}\\
&=\frac{1}{1-q}\left(\sum_{m\geq 0}\frac{q^{m^2}}{(q^2;q^2)_m}-\sum_{m\geq 0}\frac{q^{m^2+m}}{(q^2;q^2)_m}\right)\\
&=\frac{1}{1-q}\left((-q;q^2)_\infty-(-q^2;q^2)_\infty\right),
\end{align*}
where the final equation follows from \eqref{Euler-2}. The proof is complete.  \qed

\section{Proof of  Theorem \ref{thm-eo-excess}}

In this section, we aim to give a proof of  Theorem \ref{thm-eo-excess}. We find that it is equivalent to showing that
\begin{equation}\label{equiv-thm-eo-excess-1}
\sum_{n\geq 1}\left({p}^{uu}_e(n)-{p}^{uu}_o(n)\right)q^n=\frac{1}{(q^4;q^4)_\infty}\sum_{m\geq 0}\frac{q^{m^2+m}}{(q^2;q^2)_m}\frac{q^{2m+1}}{1-q^{2m+1}},
\end{equation}
\begin{equation}\label{equiv-thm-eo-excess-2}
\sum_{n\geq 1}\left({p}^{ud}_e(n)-{p}^{ud}_o(n)\right)q^n=\frac{1}{(q^2;q^2)}_\infty\sum_{m\geq 0}(-1)^mq^{2m+1}(q;q^2)_m,
\end{equation}
\begin{equation}\label{equiv-thm-eo-excess-3}
\sum_{n\geq 1}\left({p}^{du}_e(n)-{p}^{du}_o(n)\right)q^n=\sum_{m\geq 0}q^{m+1}(-q^{2m+4};q^2)_\infty,
\end{equation}
\begin{equation}\label{equiv-thm-eo-excess-4}
\sum_{n\geq 1}\left({p}^{dd}_e(n)-{p}^{dd}_o(n)\right)q^n=\sum_{m\geq 0}q^{(m+1)^2}(-q^{2m+4};q^2)_\infty.
\end{equation}

Then, we proceed show \eqref{equiv-thm-eo-excess-1}, \eqref{equiv-thm-eo-excess-2}, \eqref{equiv-thm-eo-excess-3} and \eqref{equiv-thm-eo-excess-4} in light of the smallest odd part.

{\noindent \bf Proof of \eqref{equiv-thm-eo-excess-1}.} In virtue of the smallest odd part, we get
\begin{align*}
\sum_{n\geq 1}\left({p}^{uu}_e(n)-{p}^{uu}_o(n)\right)q^n&=\sum_{n\geq 0}\frac{q^{2n+1}}{(q^{2n+1};q^2)_\infty(q^2;q^2)_n(-q^{2n+2};q^2)_\infty}\\
&=\frac{q}{(q;q^2)_\infty(-q^2;q^2)_\infty}\sum_{n\geq 0}\frac{(-q^2;q^2)_n(q;q^2)_n}{(q^2;q^2)_n(0;q^2)_n}q^{2n}\\
&=\frac{q}{(q;q^2)_\infty(-q^2;q^2)_\infty}\frac{(q;q^2)_\infty(-q^4;q^2)_\infty}{(0;q^2)_\infty(q^2;q^2)_\infty} \sum_{n\geq 0}\frac{(0;q^2)_n(q^2;q^2)_n}{(q^2;q^2)_n(-q^4;q^2)_n}q^{n}\\
&=\frac{q}{(q^4;q^4)_\infty}\sum_{n\geq 0}q^n(-q^{2n+4};q^2)_\infty\\
&=\frac{q}{(q^4;q^4)_\infty}\sum_{n\geq 0}q^n\sum_{m\geq 0}\frac{q^{(2n+4)m+2\binom{m}{2}}}{(q^2;q^2)_m}\\
&=\frac{1}{(q^4;q^4)_\infty}\sum_{m\geq 0}\frac{q^{m^2+m}}{(q^2;q^2)_m}\sum_{n\geq 0}q^{(2m+1)(n+1)}\\
&=\frac{1}{(q^4;q^4)_\infty}\sum_{m\geq 0}\frac{q^{m^2+m}}{(q^2;q^2)_m}\frac{q^{2m+1}}{1-q^{2m+1}},
\end{align*}
where the third equation follows from \eqref{Heine} and the last third equation follows from \eqref{Euler-2}. This completes the proof.  \qed

{\noindent \bf Proof of \eqref{equiv-thm-eo-excess-2}.} By considering the smallest odd part and using \eqref{q-bin-eqn}, we have
\begin{align*}
\sum_{n\geq 1}\left({p}^{ud}_e(n)-{p}^{ud}_o(n)\right)q^n&=\sum_{n\geq 0}q^{2n+1}(-q^{2n+3};q^2)_\infty\frac{1}{(q^2;q^2)_n(-q^{2n+2};q^2)_\infty}\\
&=\sum_{n\geq 0}\frac{q^{2n+1}}{(q^2;q^2)_n}\sum_{m\geq 0}\frac{(q;q^2)_m}{(q^2;q^2)_m}(-q^{2n+2})^m\\
&=\sum_{m\geq 0}\frac{(q;q^2)_m}{(q^2;q^2)_m}(-1)^mq^{2m+1}\sum_{n\geq 0}\frac{q^{(2m+2)n}}{(q^2;q^2)_n}\\
&=\sum_{m\geq 0}\frac{(q;q^2)_m}{(q^2;q^2)_m}(-1)^mq^{2m+1}\frac{1}{(q^{2m+2};q^2)_\infty}\\
&=\frac{1}{(q^2;q^2)}_\infty\sum_{m\geq 0}(-1)^mq^{2m+1}(q;q^2)_m.
\end{align*}
We arrive at \eqref{equiv-thm-eo-excess-2}, and thus the proof is complete.  \qed

{\noindent \bf Proof of \eqref{equiv-thm-eo-excess-3}.} In terms of the smallest odd part and \eqref{Heine}, we have
\begin{align*}
\sum_{n\geq 1}\left({p}^{du}_e(n)-{p}^{du}_o(n)\right)q^n&=\sum_{n\geq 0}\frac{q^{2n+1}}{(q^{2n+1};q^2)_\infty}(-q^2;q^2)_n(q^{2n+2};q^2)_\infty\\
&=q\frac{(q^2;q^2)_\infty}{(q;q^2)_\infty}\sum_{n\geq 0}\frac{(-q^2;q^2)_n(q;q^2)_n}{(q^2;q^2)_n(0;q^2)_n}q^{2n}\\
&=q\frac{(q^2;q^2)_\infty}{(q;q^2)_\infty}\frac{(q;q^2)_\infty(-q^4;q^2)_\infty}{(0;q^2)_\infty(q^2;q^2)_\infty} \sum_{m\geq 0}\frac{(0;q^2)_m(q^2;q^2)_m}{(q^2;q^2)_m(-q^4;q^2)_m}q^{m}\\
&=\sum_{m\geq 0}q^{m+1}(-q^{2m+4};q^2)_\infty.
\end{align*}
So,  \eqref{equiv-thm-eo-excess-3} is valid.  \qed

{\noindent \bf Proof of \eqref{equiv-thm-eo-excess-4}.} In view of the smallest odd part, we get
\begin{align*}
\sum_{n\geq 1}\left({p}^{dd}_e(n)-{p}^{dd}_o(n)\right)q^n&=\sum_{n\geq 0}q^{2n+1}(-q^{2n+3};q^2)_\infty(-q^2;q^2)_n(q^{2n+2};q^2)_\infty\\
&=(q^2;q^2)_\infty\sum_{n\geq 0}\frac{(-q^2;q^2)_n}{(q^2;q^2)_n}q^{2n+1}\sum_{m\geq 0}\frac{q^{(2n+3)m+2\binom{m}{2}}}{(q^2;q^2)_m}\\
&=(q^2;q^2)_\infty\sum_{m\geq 0}\frac{q^{(m+1)^2}}{(q^2;q^2)_m}\sum_{n\geq 0}\frac{(-q^2;q^2)_n}{(q^2;q^2)_n}q^{(2m+2)n}\\
&=(q^2;q^2)_\infty\sum_{m\geq 0}\frac{q^{(m+1)^2}}{(q^2;q^2)_m}\frac{(-q^{2m+4};q^2)_\infty}{(q^{2m+2};q^2)_\infty}\\
&=\sum_{m\geq 0}q^{(m+1)^2}(-q^{2m+4};q^2)_\infty,
\end{align*}
where the second equation follows from \eqref{Euler-2} and the last second equation follows from \eqref{q-bin-eqn}. So, we have proved that \eqref{equiv-thm-eo-excess-4} holds.   \qed

\section{Proofs of Theorem \ref{remin-thm-over} and Corollary \ref{remin-core-over}}

In this section, we will give the analytic proof of Theorem \ref{remin-thm-over} in Section 5.1, the combinatorial proof of Theorem \ref{remin-thm-over} in Section 5.2, and the proof of Corollary \ref{remin-core-over} in Section 5.3.

\subsection{Analytic Proof of Theorem \ref{remin-thm-over}}
In this subsection, we will give the analytic proofs of \eqref{lem-remain-excess-1} and \eqref{EE-OE-remain-eqn-1}
by considering the smallest overlined part and give the analytic proofs of \eqref{remian-eqn-3} and \eqref{remian-eqn-4} by considering the smallest non-overlined part.

{\noindent \bf Analytic proof of \eqref{lem-remain-excess-1}.} Clearly, we have
\begin{equation}\label{proof-lem-excess-1}
\sum_{n\geq 1}\left(p^{o}(n)-p^{e}(n)\right)q^n=\sum_{n\geq1}\frac{1}{(q;q)_{n-1}}\frac{q^n}{1+q^n}.
\end{equation}

In virtue of the smallest overlined part, we can get
\begin{align*}
\sum_{n\geq 1}\left(\overline{A}_{N\geq O}(n)-\overline{B}_{N\geq O}(n)\right)q^n&=\sum_{n\geq1} q^n(-q^{n+1};q)_\infty\frac{1}{(q;q)_{n-1}(-q^n;q)_\infty}\nonumber\\
&=\sum_{n\geq1}\frac{1}{(q;q)_{n-1}}\frac{q^n}{1+q^n}.
\end{align*}
Combining with \eqref{proof-lem-excess-1}, we complete the proof.  \qed

{\noindent \bf Analytic proof of \eqref{EE-OE-remain-eqn-1}.} In view of the smallest overlined part, we can get
\begin{align*}
\sum_{n\geq 1}\left(\overline{A}_{N>O}(n)-\overline{B}_{N>O}(n)\right)q^n&=\sum_{n\geq1} q^n(-q^{n+1};q)_\infty\frac{1}{(q;q)_{n}(-q^{n+1};q)_\infty}=\sum_{n\geq 1}\frac{q^n}{(q;q)_{n}}.
\end{align*}
Combining with \eqref{gen-p+}, we complete the proof.  \qed

{\noindent \bf Analytic proof of \eqref{remian-eqn-3}.} In consideration of the smallest non-overlined part, we can get
\begin{align*}
\sum_{n\geq 0}\left(\overline{A}_{O\geq N}(n)-\overline{B}_{O\geq N}(n)\right)q^n&=(q;q)_\infty+\sum_{n\geq 1} \frac{q^n}{(q^n;q)_\infty}(-q;q)_{n-1}(q^n;q)_\infty\\
&=(q;q)_\infty+\sum_{n\geq1}q^n(-q;q)_{n-1}\\
&=1+2\sum_{n\geq 1}D_e(n)q^n,
\end{align*}
where the final equation follows from \eqref{gen-d+} and \eqref{gen-d-}. The proof is complete.  \qed

{\noindent \bf Analytic proof of \eqref{remian-eqn-4}.} Clearly, we have
\begin{equation}\label{proof-remian-eqn-4}
2\sum_{n\geq 1}{H}^{o}_{ON}(n)q^n=\sum_{n\geq1}\frac{q^n}{1-q^n}\left((-q;q)_n-(q;q)_n\right).
\end{equation}
In light of the smallest non-overlined part, we can get
\begin{align*}
\sum_{n\geq0}\left(\overline{A}_{O>N}(n)-\overline{B}_{O>N}(n)\right)q^n&=(q;q)_\infty+\sum_{n\geq1}\frac{q^n}{(q^n;q)_\infty}(-q;q)_n(q^{n+1};q)_\infty\\
&=(q;q)_\infty+\sum_{n\geq1}\frac{q^n}{1-q^n}(-q;q)_n.
\end{align*}
Using \eqref{gen-d-}, we have
\begin{align*}
\sum_{n\geq0}\left(\overline{A}_{O>N}(n)-\overline{B}_{O>N}(n)\right)q^n&=1-\sum_{n\geq 1}q^n(q;q)_{n-1}+\sum_{n\geq1}\frac{q^n}{1-q^n}(-q;q)_n\\
&=1+\sum_{n\geq1}\frac{q^n}{1-q^n}\left((-q;q)_n-(q;q)_n\right).
\end{align*}
Combining with \eqref{proof-remian-eqn-4}, we complete the proof.  \qed

\subsection{Combinatorial proof of Theorem \ref{remin-thm-over}}

In this subsection, we will give the combinatorial proof of Theorem \ref{remin-thm-over}.
With the arguments in the combinatorial proofs of \eqref{lem-excess-old-1} and \eqref{EE-OE-eqn-1} in \cite[Sections 4.4 and 4.5]{Chen-He-Hu-Xie-2024}, we can give the combinatorial proofs \eqref{lem-remain-excess-1} and \eqref{EE-OE-remain-eqn-1}. The proofs are omitted here. It remains to give the combinatorial proofs of \eqref{remian-eqn-3} and \eqref{remian-eqn-4}. In  \cite[(5.7) and (5.10)]{Chen-He-Hu-Xie-2024}, Chen, He, Hu and Xie obtained that for $n\geq 1,$
\begin{equation*}\label{proof-remain-4}
	\sum_{\pi\in\overline{\mathcal{PN}}(n)}(-1)^{\ell_{O\geq N}(\pi)}=D(n),
\end{equation*}
and
\begin{equation*}\label{proof-remain-4-1}
	\sum_{\pi\in\overline{\mathcal{PN}}(n)}(-1)^{\ell_{O> N}(\pi)}={H}'_{ON}(n).
\end{equation*}

We find that in order to prove \eqref{remian-eqn-3} and \eqref{remian-eqn-4}, it suffices to show that   for $n\geq 1$,
\begin{equation}\label{proof-remain-4-000-new}\sum_{\pi\in\overline{\mathcal{D}}(n)}(-1)^{\ell_{O\geq N}(\pi)}=D_e(n)-D_o(n),
\end{equation}
and
\begin{equation}\label{proof-remin-44-05-new}
\sum_{\pi\in\overline{\mathcal{D}}(n)}(-1)^{\ell_{O>N}(\pi)}={H}^{o}_{ON}(n)-{H}^{e}_{ON}(n),
\end{equation}
where ${H}^{e}_{ON}(n)$ is the number of overpartitions $\pi$ of $n$ such that $LN(\pi)=SN(\pi)\geq 1$, $SN(\pi)\geq LO(\pi)$ and there is  an even number of overlined parts.

{\noindent \bf Combinatorial proof of \eqref{proof-remain-4-000-new}.} For $n\geq 1$, let $\pi$ be an overpartition in $\overline{\mathcal{D}}(n)$. It is clear that $\ell_{O\geq N}(\pi)=\ell(\pi)$. If we change the overlined parts in $\pi$ to non-overlined parts, then we get a distinct partition in ${\mathcal{D}}(n)$, and vice versa. So, we get
\begin{equation}\label{proof-remain-4-000}\sum_{\pi\in\overline{\mathcal{D}}(n)}(-1)^{\ell_{O\geq N}(\pi)}=\sum_{\pi\in\overline{\mathcal{D}}(n)}(-1)^{\ell(\pi)}=\sum_{\pi\in{\mathcal{D}}(n)}(-1)^{\ell(\pi)}=D_e(n)-D_o(n).
\end{equation}
 This completes the  proof. \qed

{\noindent \bf Combinatorial proof of \eqref{proof-remin-44-05-new}.} For $n\geq 1$, let $\pi$ be an overpartition in $\overline{\mathcal{D}}(n)$. It is clear that $\ell_{O>N}(\pi)=\ell_o(\pi)$. So, we get
\begin{equation}\label{proof-remin-44-05}
\sum_{\pi\in\overline{\mathcal{D}}(n)}(-1)^{\ell_{O>N}(\pi)}=\sum_{\pi\in\overline{\mathcal{D}}(n)}(-1)^{\ell_o(\pi)}.
\end{equation}

Let $\mathcal{H}'_{ON}(n)$ be the set of overpartitions counted by ${H}'_{ON}(n)$, that is, $\mathcal{H}'_{ON}(n)$ is the set of overpartitions $\pi$ of $n$ with $LN(\pi)=SN(\pi)\geq 1$ and $SN(\pi)\geq LO(\pi)$. We define an involution on $\mathcal{H}'_{ON}(n)\bigcup\overline{\mathcal{D}}(n)$ as follows. For an overpartition $\pi$ in  $\mathcal{H}'_{ON}(n)\bigcup\overline{\mathcal{D}}(n)$, we change $\pi_1$ to a non-overlined part (resp. an overlined part) if
$\pi_1$ is overlined (resp. non-overlined). This implies that
\begin{equation*}\label{proof-remin-44-02}
\sum_{\pi\in\mathcal{H}'_{ON}(n)\bigcup\overline{\mathcal{D}}(n)}(-1)^{\ell_o(\pi)}=0.
\end{equation*}
Combining with \eqref{proof-remin-44-05}, we have
\begin{equation*}\label{proof-remin-44-04}
\sum_{\pi\in\overline{\mathcal{D}}(n)}(-1)^{\ell_{O>N}(\pi)}=-\sum_{\pi\in\mathcal{H}'_{ON}(n)}(-1)^{\ell_o(\pi)}={H}^{o}_{ON}(n)-{H}^{e}_{ON}(n).
\end{equation*}

The proof is complete.  \qed

\subsection{Proof of Corollary \ref{remin-core-over}}

The objective of this subsection is to give the proof of Corollary \ref{remin-core-over}. Clearly, \eqref{core-remain-excess-1}, \eqref{core-remain-excess-2} and \eqref{core-remain-excess-3} immediately follows from \eqref{EE-OE-remain-eqn-1}, \eqref{remian-eqn-3} and \eqref{remian-eqn-4} respectively. It remains to show \eqref{core-remain-excess-4}. Appealing to \eqref{lem-remain-excess-1}, we find that it is equivalent to showing that
for $n\geq 1$,
\begin{equation}\label{lem-remain-cor-excess-1}
p^{o}(n)-p^{e}(n)\geq 0\text{ with strict inequality if }n\neq 2.
\end{equation}

For $n\geq 1$, let $\hat{p}^o(n)$ be the number of partitions $\pi$ of $n$ such that  $\pi_1$ is even, $\pi_2\leq \frac{\pi_1}{2}$ and $\frac{\pi_1}{2}$ appears an even number of times. For a partition $\pi$ counted by  $\hat{p}^o(n)$, it is clear that the largest part of $\pi$ appears once in $\pi$, and so $\pi$ is a partition enumerated by ${p}^o(n)$.
In order to show \eqref{lem-remain-cor-excess-1}, it suffices to prove that
\begin{itemize}
\item[(1)] $\hat{p}^o(n)=p^{e}(n)$;

\item[(2)] there exists a partition counted by  ${p}^o(n)$ but not enumerated by $\hat{p}^o(n)$ for $n\neq 2$.
\end{itemize}

We will give an analytic proof and a combinatorial proof of the condition (1), and then we will give a proof of the condition (2).

{\noindent \bf Analytic proof of the condition (1).} By considering the largest part, we get
\[\sum_{n\geq 1}\hat{p}^o(n)q^n=\sum_{n\geq 1}\frac{q^{2n}}{1-q^{2n}}\frac{1}{(q;q)_{n-1}}=\sum_{n\geq 1}{p}^e(n)q^n.\]
This completes the proof.   \qed

{\noindent \bf Combinatorial proof of the condition (1).}   For $n\geq 1$, let $\pi$ be a partition counted by $\hat{p}^o(n)$. We set
\[\lambda=\left(\frac{\pi_1}{2},\frac{\pi_1}{2},\pi_2,\ldots,\pi_{\ell(\pi)}\right).\]
Clearly, $\lambda$ is a partition enumerated by ${p}^e(n)$.

Conversely, for a partition $\lambda$ enumerated by ${p}^e(n)$, we set
\[\pi=\left(2\lambda_1,\lambda_3,\ldots,\lambda_{\ell(\lambda)}\right).\]
Obviously, $\pi$ is a partition counted by  $\hat{p}^o(n)$. This completes the proof.  \qed

{\noindent \bf Proof of the condition (2).} For $n\geq 1$, let $\pi$ be a partition counted by $\hat{p}^o(n)$. By definition, we know that the largest part of $\pi$ is even. So, we just need to find a partition $\lambda$ counted by ${p}^o(n)$ with the largest part being odd for $n\neq 2$. There are two cases.

Case 1: $n$ is odd. In such case, we set $\lambda=(\underbrace{1,\ldots,1}_{n's\ 1})$.

Case 2: $n$ is even and $n\geq 4$. In such case, we set $\lambda=(3,\underbrace{1,\ldots,1}_{(n-3)'s\ 1})$.

 In either case, we find a partition $\lambda$  counted by  ${p}^o(n)$ but not enumerated by $\hat{p}^o(n)$, and thus the proof is complete.  \qed

\section{Analytic proof of Theorem \ref{new-A-B-thm}}

This section is devoted to giving the analytic proof of Theorem \ref{new-A-B-thm}. We will give the analytic proofs of \eqref{E-O-EQN(1)}-\eqref{E-O-EQN(4)}
by considering the largest non-overlined part and give the analytic proofs of \eqref{E-O-EQN(5)}-\eqref{E-O-EQN(8)} by considering the largest overlined part.

{\noindent \bf Analytic proofs of \eqref{E-O-EQN(1)} and \eqref{E-O-EQN(2)}.} Appealing to \eqref{gen-d-}, we have
\begin{align*}\sum_{n\geq 0}\left(A_{O\leq N}(n)-B_{O\leq N}(n)\right)q^n&=(q;q)_\infty+\sum_{n\geq 1}\left(\widetilde{A}_{O\leq N}(n)-\widetilde{B}_{O\leq N}(n)\right)q^n\\
&=\sum_{n\geq0}\left(D_e(n)-D_o(n)\right)q^{n}+\sum_{n\geq 1}\left(\widetilde{A}_{O\leq N}(n)-\widetilde{B}_{O\leq N}(n)\right)q^n.
\end{align*}

So, we just need to show that
\begin{equation}\label{proof-new-OLEQn}
\sum_{n\geq 1}\left(\widetilde{A}_{O\leq N}(n)-\widetilde{B}_{O\leq N}(n)\right)q^n=\sum_{n\geq 1}D(n)q^n.
\end{equation}

In virtue of the largest non-overlined part, we can get
\begin{equation*}
\sum_{n\geq 1}\left(\widetilde{A}_{O\leq N}(n)-\widetilde{B}_{O\leq N}(n)\right)q^n=\sum_{n\geq 1}\frac{q^n}{(q;q)_n}(-q^{n+1};q)_\infty(q;q)_n=\sum_{n\geq 1}q^n(-q^{n+1};q)_\infty.
\end{equation*}
Combining with \eqref{gen-d+}, we arrive at \eqref{proof-new-OLEQn}.  The proof is complete. \qed

{\noindent \bf Analytic proofs of \eqref{E-O-EQN(3)} and \eqref{E-O-EQN(4)}.} Clearly, we have
\begin{equation}\label{proof-E-O-EQN(4)}
\sum_{n\geq 1}H_{O< N}(n)q^n=\sum_{n\geq 1}\frac{q^n}{1-q^n}(-q^n;q)_\infty,
\end{equation}
and
\begin{equation}\label{proof-E-O-EQN(3)}
2\sum_{n\geq 1}H_{O< N}^o(n)q^n=\sum_{n\geq 1}\frac{q^n}{1-q^n}\left((-q^n;q)_\infty-(q^n;q)_\infty\right).
\end{equation}

 In view of the largest non-overlined part, we can get
\[
\sum_{n\geq 1}\left(\widetilde{A}_{O<N}(n)-\widetilde{B}_{O<N}(n)\right)q^n=\sum_{n\geq 1}\frac{q^n}{(q;q)_n}(-q^{n};q)_\infty(q;q)_{n-1}
=\sum_{n\geq 1}\frac{q^n}{1-q^n}(-q^n;q)_\infty.
\]
Combining with \eqref{proof-E-O-EQN(4)}, we prove that \eqref{E-O-EQN(4)} holds. Using \eqref{gen-d-}, we get
\begin{align*}
\sum_{n\geq 0}\left(A_{O<N}(n)-B_{O<N}(n)\right)q^n&=(q;q)_\infty+\sum_{n\geq 1}\left(\widetilde{A}_{O<N}(n)-\widetilde{B}_{O<N}(n)\right)q^n\\
&=(q;q)_\infty+\sum_{n\geq 1}\frac{q^n}{1-q^n}(-q^n;q)_\infty\\
&=1-\sum_{n\geq 1}q^n(q^{n+1};q)_\infty+\sum_{n\geq 1}\frac{q^n}{1-q^n}(-q^n;q)_\infty\\
&=1+\sum_{n\geq 1}\frac{q^n}{1-q^n}\left((-q^n;q)_\infty-(q^n;q)_\infty\right).
\end{align*}
Combining with \eqref{proof-E-O-EQN(3)}, we derive that \eqref{E-O-EQN(3)} is valid. The proof is complete. \qed

{\noindent \bf Analytic proofs of \eqref{E-O-EQN(5)} and \eqref{E-O-EQN(6)}.} Clearly, we have
\begin{equation}\label{proof-E-O-EQN(6)}
\sum_{n\geq 1}\left(p'_{o}(n)-p'_{e}(n)\right)q^n=\sum_{n\geq 1}\frac{q^n}{1+q^n}\frac{1}{(q^{n+1};q)_\infty},
\end{equation}
and
\begin{equation}\label{proof-E-O-EQN(5)}
2\sum_{n\geq 1}\left(p''_{o}(n)-p''_{e}(n)\right)q^n=\sum_{n\geq 1}\frac{q^n}{1+q^n}\left(\frac{1}{(q^{n+1};q)_\infty}-\frac{1}{(-q^{n+1};q)_\infty}\right).
\end{equation}

 In consideration of the largest overlined part, we can get
\begin{align*}
\sum_{n\geq 1}\left(\widetilde{A}_{N\leq O}(n)-\widetilde{B}_{N\leq O}(n)\right)q^n&=\sum_{n\geq 1}q^n(-q;q)_{n-1}\frac{1}{(-q;q)_n(q^{n+1};q)_\infty}\\
&=\sum_{n\geq 1}\frac{q^n}{1+q^n}\frac{1}{(q^{n+1};q)_\infty}.
\end{align*}
Combining with \eqref{proof-E-O-EQN(6)}, we prove that \eqref{E-O-EQN(6)} holds. Using \eqref{gen-p-}, we get
\begin{align*}
\sum_{n\geq 0}\left(A_{N\leq O}(n)-B_{N\leq O}(n)\right)q^n&=\frac{1}{(-q;q)_\infty}+\sum_{n\geq 1}\left(\widetilde{A}_{N\leq O}(n)-\widetilde{B}_{N\leq O}(n)\right)q^n\\
&=1-\sum_{n\geq 1}\frac{q^n}{(-q^{n};q)_\infty}+\sum_{n\geq 1}\frac{q^n}{1+q^n}\frac{1}{(q^{n+1};q)_\infty}\\
&=1+\sum_{n\geq 1}\frac{q^n}{1+q^n}\left(\frac{1}{(q^{n+1};q)_\infty}-\frac{1}{(-q^{n+1};q)_\infty}\right).
\end{align*}
Combining with \eqref{proof-E-O-EQN(5)}, we derive that \eqref{E-O-EQN(5)} is valid. The proof is complete. \qed

{\noindent \bf Analytic proofs of \eqref{E-O-EQN(7)} and \eqref{E-O-EQN(8)}.} Appealing to \eqref{gen-p-}, we have
\begin{align*}\sum_{n\geq 0}\left(A_{N<O}(n)-B_{N<O}(n)\right)q^n&=\frac{1}{(-q;q)_\infty}+\sum_{n\geq 1}\left(\widetilde{A}_{N<O}(n)-\widetilde{B}_{N<O}(n)\right)q^n\\
&=\sum_{n\geq0}\left(p_e(n)-p_o(n)\right)q^{n}+\sum_{n\geq 1}\left(\widetilde{A}_{N<O}(n)-\widetilde{B}_{N<O}(n)\right)q^n.
\end{align*}

So, we just need to show that
\begin{equation}\label{proof-new-OLEQn-<}
\sum_{n\geq 1}\left(\widetilde{A}_{N<O}(n)-\widetilde{B}_{N<O}(n)\right)q^n=\sum_{n\geq 1}p(n)q^n.
\end{equation}

In light of the largest overlined part, we can get
\begin{equation*}
\sum_{n\geq 1}\left(\widetilde{A}_{N<O}(n)-\widetilde{B}_{N<O}(n)\right)q^n=\sum_{n\geq 1}q^n(-q;q)_{n-1}\frac{1}{(-q;q)_{n-1}(q^{n};q)_\infty}\\
=\sum_{n\geq 1}\frac{q^n}{(q^{n};q)_\infty}.
\end{equation*}
Combining with \eqref{gen-p+}, we arrive at \eqref{proof-new-OLEQn-<}.  The proof is complete. \qed

\section{Combinatorial proof of Theorem \ref{new-A-B-thm}}

The objective of this section is to give the combinatorial proof of Theorem \ref{new-A-B-thm}. We will give the combinatorial proofs of  \eqref{E-O-EQN(1)} and  \eqref{E-O-EQN(2)} in Section 7.1, \eqref{E-O-EQN(3)} and  \eqref{E-O-EQN(4)} in Section 7.2, and \eqref{E-O-EQN(5)}-\eqref{E-O-EQN(8)} in Section 7.3.

\subsection{Combinatorial proofs of  \eqref{E-O-EQN(1)} and  \eqref{E-O-EQN(2)}}

In this subsection, we aim to give the combinatorial proofs of \eqref{E-O-EQN(1)} and  \eqref{E-O-EQN(2)}, which are equivalent to showing that for $n\geq 1$,
\begin{equation}\label{E-O-EQN(1)-equiv}
\sum_{\pi\in\overline{\mathcal{P}}(n)}(-1)^{\ell_{O\leq N}(\pi)}=2D_e(n),
\end{equation}
and
\begin{equation}\label{E-O-EQN(2)-equiv}
\sum_{\pi\in\overline{\mathcal{PN}}(n)}(-1)^{\ell_{O\leq N}(\pi)}=D(n).
\end{equation}

To do this, we introduce the following notations.
\begin{itemize}
\item Let $\mathcal{C}_{O\leq N}(n)$ be the set of overpartitions $\pi$ of $n$ such that
there are at least one non-overlined part and at least one overlined part in $\pi$ and $\widetilde{SN}(\pi)\geq \widetilde{SO}(\pi)$.

\item Let $\mathcal{F}_{O\leq N}(n)$ be the set of overpartitions $\pi$ of $n$ such that there are at least two non-overlined parts in $\pi$ and $\widetilde{SO}(\pi)>\widetilde{SN}(\pi)$.

\item Let $\mathcal{H}_{O\leq N}(n)$ be the set of overpartitions $\pi$ of $n$ such that there is exactly one non-overlined part in $\pi$ and $\widetilde{SO}(\pi)>\widetilde{SN}(\pi)$.
\end{itemize}
Clearly, we have
\[\overline{\mathcal{PN}}(n)=\mathcal{C}_{O\leq N}(n)\bigcup\mathcal{F}_{O\leq N}(n)\bigcup\mathcal{H}_{O\leq N}(n).\]

Then, we proceed to present the combinatorial proofs of \eqref{E-O-EQN(1)-equiv} and  \eqref{E-O-EQN(2)-equiv}.

{\noindent \bf Combinatorial proofs of \eqref{E-O-EQN(1)-equiv} and  \eqref{E-O-EQN(2)-equiv}.}
By restricting the involution $\varphi$ defined in Definition \ref{defi-involution-varphi} on $\mathcal{C}_{O\leq N}(n)\bigcup\mathcal{F}_{O\leq N}(n)$, we get
\begin{equation}\label{proof-E-O-EQN(1)-equiv-CF}
 \sum_{\pi\in\mathcal{C}_{O\leq N}(n)\bigcup\mathcal{F}_{O\leq N}(n)}(-1)^{\ell_{O\leq N}(\pi)}=0.
\end{equation}

Let $\pi=(\pi_1,\pi_2,\ldots,\pi_m)$ be an overpartition in $\mathcal{H}_{O\leq N}(n)$. Then, we have  $\widetilde{LN}(\pi)=\widetilde{SN}(\pi)<\widetilde{SO}(\pi)$, which yields $\ell_{O\leq N}(\pi)=0$, and so $(-1)^{\ell_{O\leq N}(\pi)}=1$. If we change the overlined parts $\pi_1,\ldots,\pi_{m-1}$ in $\pi$ to non-overlined parts, then we get a partition in $\mathcal{D}(n)$, and vice versa. This implies that the number of overpartitions in $\mathcal{H}_{O\leq N}(n)$ is $D(n)$. So, we get
\begin{equation}\label{proof-E-O-EQN(1)-equiv-H-2}
 \sum_{\pi\in\mathcal{H}_{O\leq N}(n)}(-1)^{\ell_{O\leq N}(\pi)}=\sum_{\pi\in\mathcal{H}_{O\leq N}(n)}1=D(n).
\end{equation}
Combining with \eqref{proof-E-O-EQN(1)-equiv-CF}, we arrive at \eqref{E-O-EQN(2)-equiv}.

For an overpartition $\pi\in \overline{\mathcal{D}}(n)$, it is clear that $\ell_{O\leq N}(\pi)=\ell(\pi)$. Using \eqref{proof-remain-4-000}, we get
\begin{equation*}\label{proof-E-O-EQN(1)-equiv-D}\sum_{\pi\in\overline{\mathcal{D}}(n)}(-1)^{\ell_{O\leq N}(\pi)}=\sum_{\pi\in\overline{\mathcal{D}}(n)}(-1)^{\ell(\pi)}=D_e(n)-D_o(n).
\end{equation*}
Combining with \eqref{proof-E-O-EQN(1)-equiv-CF} and \eqref{proof-E-O-EQN(1)-equiv-H-2}, we deduce that \eqref{E-O-EQN(1)-equiv} is valid.  The proof is complete.  \qed

\subsection{Combinatorial proofs of  \eqref{E-O-EQN(3)} and  \eqref{E-O-EQN(4)}}

In this subsection, we aim to give the combinatorial proofs of \eqref{E-O-EQN(3)} and  \eqref{E-O-EQN(4)}, which are equivalent to showing that for $n\geq 1$,
\begin{equation}\label{E-O-EQN(3)-equiv}
\sum_{\pi\in\overline{\mathcal{P}}(n)}(-1)^{\ell_{O< N}(\pi)}=2H_{O<N}^o(n),
\end{equation}
and
\begin{equation}\label{E-O-EQN(4)-equiv}
\sum_{\pi\in\overline{\mathcal{PN}}(n)}(-1)^{\ell_{O< N}(\pi)}=H_{O<N}(n).
\end{equation}

To do this, we introduce the following notations.
\begin{itemize}
\item Let $\mathcal{C}_{O< N}(n)$ be the set of overpartitions $\pi$ of $n$ such that
there are at least one non-overlined part and at least one overlined part in $\pi$, $\widetilde{SN}(\pi)\geq \widetilde{SO}(\pi)$ and $\widetilde{LN}(\pi)>\widetilde{SO}(\pi)$.

\item Let $\mathcal{F}_{O< N}(n)$ be the set of overpartitions $\pi$ of $n$ such that there are at least two non-overlined parts in $\pi$, $\widetilde{SO}(\pi)>\widetilde{SN}(\pi)$ and $\widetilde{LN}(\pi)> \widetilde{SN}(\pi)$.

\item Let $\mathcal{H}_{O< N}(n)$ be the set of overpartitions $\pi$ of $n$ such that there is at least one non-overlined part in $\pi$ and $\widetilde{SO}(\pi)\geq\widetilde{LN}(\pi)=\widetilde{SN}(\pi)$, that is, $\mathcal{H}_{O< N}(n)$ is the set of overpartitions counted by ${H}_{O< N}(n)$.
\end{itemize}
Clearly, we have
\[\overline{\mathcal{PN}}(n)=\mathcal{C}_{O< N}(n)\bigcup\mathcal{F}_{O< N}(n)\bigcup\mathcal{H}_{O< N}(n).\]

Then, we proceed to present the combinatorial proofs of \eqref{E-O-EQN(3)-equiv} and  \eqref{E-O-EQN(4)-equiv}.

{\noindent \bf Combinatorial proofs of \eqref{E-O-EQN(3)-equiv} and  \eqref{E-O-EQN(4)-equiv}.}
By restricting the involution $\varphi$ defined in Definition \ref{defi-involution-varphi} on $\mathcal{C}_{O< N}(n)\bigcup\mathcal{F}_{O< N}(n)$, we get
\begin{equation}\label{proof-E-O-EQN(3)-equiv-CF}
 \sum_{\pi\in\mathcal{C}_{O< N}(n)\bigcup\mathcal{F}_{O< N}(n)}(-1)^{\ell_{O< N}(\pi)}=0.
\end{equation}

For an overpartition $\pi\in\mathcal{H}_{O<N}(n)$, we have  $\widetilde{LN}(\pi)\leq\widetilde{SO}(\pi)$, which yields $\ell_{O< N}(\pi)=0$, and so $(-1)^{\ell_{O< N}(\pi)}=1$. Then, we have
\begin{equation}\label{proof-E-O-EQN(3)-equiv-H}
 \sum_{\pi\in\mathcal{H}_{O< N}(n)}(-1)^{\ell_{O<N}(\pi)}=\sum_{\pi\in\mathcal{H}_{O< N}(n)}1={H}_{O< N}(n).
\end{equation}
Combining with \eqref{proof-E-O-EQN(3)-equiv-CF}, we arrive at \eqref{E-O-EQN(4)-equiv}.

By \eqref{proof-E-O-EQN(3)-equiv-CF} and \eqref{proof-E-O-EQN(3)-equiv-H}, we find that in order to show \eqref{E-O-EQN(3)-equiv}, it remains to prove that
\begin{equation}\label{proof-E-O-EQN(3)-equiv-D}
\sum_{\pi\in\overline{\mathcal{D}}(n)}(-1)^{\ell_{O<N}(\pi)}=-\sum_{\pi\in\mathcal{H}_{O< N}(n)}(-1)^{\ell_{o}(\pi)}.
\end{equation}

For an overpartition $\pi\in \overline{\mathcal{D}}(n)$, it is clear that $\ell_{O<N}(\pi)=\ell_o(\pi)$, and so
\begin{equation}\label{O<N-2-EQN-1}
\sum_{\pi\in\overline{\mathcal{D}}(n)}(-1)^{\ell_{O<N}(\pi)}=\sum_{\pi\in\overline{\mathcal{D}}(n)}(-1)^{\ell_{o}(\pi)}.
\end{equation}

Again by restricting the involution $\varphi$ defined in Definition \ref{defi-involution-varphi} on $\mathcal{H}_{O< N}(n)\bigcup\overline{\mathcal{D}}(n)$,
we get
\begin{equation*}
\sum_{\pi\in\mathcal{H}_{O< N}(n)\bigcup\overline{\mathcal{D}}(n)}(-1)^{\ell_{o}(\pi)}=0.
\end{equation*}
Combining with \eqref{O<N-2-EQN-1}, we arrive at \eqref{proof-E-O-EQN(3)-equiv-D}, and thus the proof is complete. \qed

\subsection{Combinatorial proofs of  \eqref{E-O-EQN(5)}-\eqref{E-O-EQN(8)}}

In this subsection, we aim to give the combinatorial proofs of \eqref{E-O-EQN(5)}-\eqref{E-O-EQN(8)}, which are equivalent to showing that for $n\geq 1$,
\begin{equation}\label{E-O-EQN(5)-equiv}
\sum_{\pi\in\overline{\mathcal{P}}(n)}(-1)^{\ell_{N\leq O}(\pi)}=2\left(p''_{o}(n)-p''_{e}(n)\right),
\end{equation}
\begin{equation}\label{E-O-EQN(6)-equiv}
\sum_{\pi\in\overline{\mathcal{PO}}(n)}(-1)^{\ell_{N\leq O}(\pi)}=p'_{o}(n)-p'_{e}(n),
\end{equation}
\begin{equation}\label{E-O-EQN(7)-equiv}
\sum_{\pi\in\overline{\mathcal{P}}(n)}(-1)^{\ell_{N<O}(\pi)}=2p_e(n),
\end{equation}
\begin{equation}\label{E-O-EQN(8)-equiv}
\sum_{\pi\in\overline{\mathcal{PO}}(n)}(-1)^{\ell_{N<O}(\pi)}=p(n).
\end{equation}

To do this, we introduce the following notations.
\begin{itemize}
\item Let $\widetilde{\mathcal{C}}_{NO}(n)$ be the set of overpartitions $\pi$ of $n$ such that
there are at least one non-overlined part and at least one overlined part in $\pi$ and $\widetilde{SN}(\pi)<\widetilde{SO}(\pi)$.

\item Let $\widetilde{\mathcal{F}}_{NO}(n)$ be the set of overpartitions $\pi$ of $n$ such that there are at least two overlined parts in $\pi$ and $\widetilde{SO}(\pi)\leq\widetilde{SN}(\pi)$.

\item Let $\widetilde{\mathcal{H}}_{NO}(n)$ be the set of overpartitions $\pi$ of $n$ such that there is exactly one overlined part in $\pi$ and $\widetilde{SO}(\pi)\leq\widetilde{SN}(\pi)$.
\end{itemize}
Clearly, we have
\[\overline{\mathcal{PO}}(n)=\widetilde{\mathcal{C}}_{NO}(n)\bigcup\widetilde{\mathcal{F}}_{NO}(n)\bigcup\widetilde{\mathcal{H}}_{NO}(n).\]

By restricting the involution $\varphi$ defined in Definition \ref{defi-involution-varphi} on $\widetilde{\mathcal{C}}_{NO}(n)\bigcup\widetilde{\mathcal{F}}_{NO}(n)$,  we get
\begin{equation}\label{proof-E-O-EQN(5-8)-equiv-CF-leq}
 \sum_{\pi\in\widetilde{\mathcal{C}}_{NO}(n)\bigcup\widetilde{\mathcal{F}}_{NO}(n)}(-1)^{\ell_{N\leq O}(\pi)}=0,
\end{equation}
and
\begin{equation}\label{proof-E-O-EQN(5-8)-equiv-CF-<}
 \sum_{\pi\in\widetilde{\mathcal{C}}_{NO}(n)\bigcup\widetilde{\mathcal{F}}_{NO}(n)}(-1)^{\ell_{N<O}(\pi)}=0.
\end{equation}

For a partition $\pi\in {\mathcal{P}}(n)$, it is clear that $\ell_{N\leq O}(\pi)=\ell_{N<O}(\pi)=\ell(\pi)$, and so
\begin{equation*}
\sum_{\pi\in{\mathcal{P}}(n)}(-1)^{\ell_{N\leq O}(\pi)}=\sum_{\pi\in{\mathcal{P}}(n)}(-1)^{\ell_{N<O}(\pi)}=\sum_{\pi\in{\mathcal{P}}(n)}(-1)^{\ell(\pi)}=p_e(n)-p_o(n).
\end{equation*}
Combining with \eqref{proof-E-O-EQN(5-8)-equiv-CF-leq} and \eqref{proof-E-O-EQN(5-8)-equiv-CF-<}, we find that in order to show \eqref{E-O-EQN(5)-equiv}-\eqref{E-O-EQN(8)-equiv}, it remains to prove that for $n\geq 1$,
\begin{equation}\label{new-proof-6-PH-remain}
\sum_{\pi\in\widetilde{\mathcal{H}}_{NO}(n)}(-1)^{\ell_{N\leq O}(\pi)}=p'_{o}(n)-p'_{e}(n),
\end{equation}
and
\begin{equation}\label{new-proof-78-H-1}
 \sum_{\pi\in\widetilde{\mathcal{H}}_{NO}(n)}(-1)^{\ell_{N<O}(\pi)}=p(n).
\end{equation}

Then, we proceed to present the combinatorial proofs of \eqref{new-proof-6-PH-remain} and \eqref{new-proof-78-H-1}.

{\noindent \bf Combinatorial proof of \eqref{new-proof-6-PH-remain}.} For $n\geq 1$, let $\pi$ be an overpartition in $\widetilde{\mathcal{H}}_{NO}(n)$.
It is clear that $\ell_{N\leq O}(\pi)$ is the number of non-overlined parts of size $\widetilde{SO}(\pi)$  in $\pi$. This implies that
$(-1)^{\ell_{N\leq O}(\pi)}=1$ (resp. $(-1)^{\ell_{N\leq O}(\pi)}=-1$) if there is an even (resp. odd) number of non-overlined parts of size $\widetilde{SO}(\pi)$ in $\pi$.

By restricting the involution $\varphi$ defined in Definition \ref{defi-involution-varphi} on $\widetilde{\mathcal{H}}_{NO}(n)\bigcup\mathcal{P}(n)$, we know that the map $\varphi$ is a bijection between $\widetilde{\mathcal{H}}_{NO}(n)$ and $\mathcal{P}(n)$. Furthermore, we find that the number of overpartitions $\pi$ in $\widetilde{\mathcal{H}}_{NO}(n)$ with an even (resp. odd) number of non-overlined parts of size $\widetilde{SO}(\pi)$ equals the number of partitions in $\mathcal{P}(n)$ such that the smallest part appears an odd (resp. even) number of times. This completes the proof.  \qed

{\noindent \bf Combinatorial proof of \eqref{new-proof-78-H-1}.}  For $n\geq 1$, let $\pi$ be an overpartition in $\widetilde{\mathcal{H}}_{NO}(n)$. Then, we have  $\widetilde{LO}(\pi)=\widetilde{SO}(\pi)\leq\widetilde{SN}(\pi)$, which yields $\ell_{N< O}(\pi)=0$, and so $(-1)^{\ell_{N< O}(\pi)}=1$. Recalling that the map $\varphi$ defined in Definition \ref{defi-involution-varphi} is a bijection between $\widetilde{\mathcal{H}}_{NO}(n)$ and $\mathcal{P}(n)$, we get
\begin{equation*}
 \sum_{\pi\in\widetilde{\mathcal{H}}_{NO}(n)}(-1)^{\ell_{N< O}(\pi)}=\sum_{\pi\in\widetilde{\mathcal{H}}_{NO}(n)}1=p(n).
\end{equation*}

 The proof is complete.  \qed

\section{Combinatorial proof of Corollary \ref{coro-equal-not}}

This section is devoted to giving  the combinatorial proof of Corollary \ref{coro-equal-not}. Clearly, it suffices to show the following lemma.
\begin{lem}
For $n\geq 1$, there is an involution $\phi$ on $\overline{\mathcal{P}}(n)$ and $\phi$ is also an involution on $\overline{\mathcal{PN}}(n)$. Moreover, for an overpartition $\pi\in\overline{\mathcal{P}}(n)$, let $\lambda=\phi(\pi)$. We have
\[\ell_{O\leq N}(\lambda)\equiv \ell_{O\geq N}(\pi)\pmod 2.\]
\end{lem}

\pf Let $\pi$ be an overpartition in $\overline{\mathcal{P}}(n)$. We define  $\phi\colon \pi \rightarrow \lambda$
   as follows. There are two cases.

 \begin{itemize}
 \item Case 1: $\ell_{O\leq N}(\pi)\equiv \ell_{O\geq N}(\pi)\pmod 2$. In this case, we set $\lambda=\pi$. Clearly, we have $\lambda\in\overline{\mathcal{P}}(n)$ and  $\ell_{O\leq N}(\lambda)=\ell_{O\leq N}(\pi)\equiv \ell_{O\geq N}(\pi)\pmod 2$.

 \item Case 2: $\ell_{O\leq N}(\pi)\not\equiv \ell_{O\geq N}(\pi)\pmod 2$. In this case, non-overlined parts appear in $\pi$, and so $SN(\pi)=\widetilde{SN}(\pi)$ and $LN(\pi)=\widetilde{LN}(\pi)$. There are two subcases.

 \begin{itemize}
 \item  Case 2.1: there is exactly one part of size ${SN}(\pi)$ in $\pi$, there is exactly one part of size ${LN}(\pi)$ in $\pi$, and there is no part of size $t$ in $\pi$ for ${SN}(\pi)< t< {LN}(\pi)$.

     \begin{itemize}
    \item If $\ell_{O\geq N}(\pi)$ is odd, then we assume that $k$ is the smallest integer such that $k>{LN}(\pi)$ and $\overline{k}$ appears in $\pi$. The overpartition $\lambda$ is obtained by changing the smallest non-overlined part in $\pi$ to an overlined part and changing $\overline{k}$ in $\pi$ to $k$. Clearly, we have $\lambda\in\overline{\mathcal{P}}(n)$ and $\ell_{O\leq N}(\lambda)=\ell_{O\leq N}(\pi)+1\equiv \ell_{O\geq N}(\pi)\pmod 2$.

     \item If $\ell_{O\leq N}(\pi)$ is odd, then we assume that $k$ is the largest integer such that $k<{SN}(\pi)$ and $\overline{k}$ appears in $\pi$. The overpartition $\lambda$ is obtained by changing the largest non-overlined part in $\pi$ to an overlined part and changing $\overline{k}$ in $\pi$ to $k$. Clearly, we have $\lambda\in\overline{\mathcal{P}}(n)$ and $\ell_{O\leq N}(\lambda)=\ell_{O\leq N}(\pi)-1\equiv \ell_{O\geq N}(\pi)\pmod 2$.

     \end{itemize}

      \item  Case 2.2: there are at least two parts of size ${SN}(\pi)$ in $\pi$, or there are at least two parts of size ${LN}(\pi)$ in $\pi$, or there is at least one part of size greater than ${SN}(\pi)$ and less than ${LN}(\pi)$ in $\pi$. We set $k$ be the smallest integer such that $k>{SN}(\pi)$ and there is at least one part of size $k$ in $\pi$ if there is exactly one part of size ${SN}(\pi)$ in $\pi$, and set $k={SN}(\pi)$ otherwise.

     \begin{itemize}
    \item If $\overline{k}$ occurs in $\pi$, then $\lambda$ is obtained by changing $\overline{k}$ in $\pi$ to $k$. Clearly, we have $\lambda\in\overline{\mathcal{P}}(n)$ and $\ell_{O\leq N}(\lambda)=\ell_{O\leq N}(\pi)-1\equiv \ell_{O\geq N}(\pi)\pmod 2$.

         \item If $\overline{k}$ does not occur in $\pi$, then $\lambda$ is obtained by changing a $k$ in $\pi$ to $\overline{k}$. Clearly, we have $\lambda\in\overline{\mathcal{P}}(n)$ and $\ell_{O\leq N}(\lambda)=\ell_{O\leq N}(\pi)+1\equiv \ell_{O\geq N}(\pi)\pmod 2$.
            \end{itemize}

 \end{itemize}

\end{itemize}

In conclusion, we have $\lambda=\phi(\pi)\in\overline{\mathcal{P}}(n)$ and $\ell_{O\leq N}(\lambda)\equiv \ell_{O\geq N}(\pi)\pmod 2$. It can be checked that $\phi$ is an involution on $\overline{\mathcal{P}}(n)$.
Moreover, we have $\lambda=\phi(\pi)\in\overline{\mathcal{PN}}(n)$ if $\pi\in\overline{\mathcal{PN}}(n)$, which implies that $\phi$ is  an involution on $\overline{\mathcal{PN}}(n)$. The proof is complete. \qed

We present four examples to illustrate the map $\phi$ above.
\begin{itemize}

\item[(1)] For an overparition $\pi=(\overline{8},7,5,\overline{3},\overline{2})$, we have $LN(\pi)=7$, $SN(\pi)=5$, $\ell_{O\leq N}(\pi)=2$ and $\ell_{O\geq N}(\pi)=1$. Moreover, there is exactly one part of size $5$ in $\pi$, there is exactly one part of size $7$ in $\pi$, and there is no part of size $t$ in $\pi$ for $5< t< 7$. Note that $\ell_{O\geq N}(\pi)$ is odd, then we replace $5$ and $\overline{8}$ in $\pi$ by $\overline{5}$ and $8$ respectively to get $\lambda=({8},7,\overline{5},\overline{3},\overline{2})$. Clearly, we have  $\ell_{O\leq N}(\lambda)=\ell_{O\leq N}(\pi)+1=3$.

\item[(2)] For an overparition $\pi=(7,\overline{5},\overline{3},\overline{2})$, we have $LN(\pi)=SN(\pi)=7$, $\ell_{O\leq N}(\pi)=3$ and $\ell_{O\geq N}(\pi)=0$. Moreover, there is exactly one part of size $7$ in $\pi$. Note that $\ell_{O\leq N}(\pi)$ is odd, then we replace $7$ and $\overline{5}$ in $\pi$ by $\overline{7}$ and $5$ respectively to get $\lambda=(\overline{7},{5},\overline{3},\overline{2})$. Clearly, we have  $\ell_{O\leq N}(\lambda)=\ell_{O\leq N}(\pi)-1=2$.

\item[(3)] For an overparition $\pi=(\overline{13},\overline{11},\overline{10},10,10,\overline{9},8,8,7,\overline{5},\overline{3},\overline{2})$, we have $LN(\pi)=10$, $SN(\pi)=7$, $\ell_{O\leq N}(\pi)=5$ and $\ell_{O\geq N}(\pi)=4$. Moreover, there are three parts $\overline{10},10,10$ of size $10$ in $\pi$ and there exist parts $\overline{9},8,8$ of size greater than $7$ and less than $10$ in $\pi$. Note that there is exactly one part of size $7$ in $\pi$ and $\overline{8}$ does not occur in $\pi$, then we replace a ${8}$ in $\pi$ by $\overline{8}$ to get $\lambda=(\overline{13},\overline{11},\overline{10},10,10,\overline{9},\overline{8},8,7,\overline{5},\overline{3},\overline{2})$. Clearly, we have  $\ell_{O\leq N}(\lambda)=\ell_{O\leq N}(\pi)+1=6$.

\item[(4)] For an overparition $\pi=(\overline{13},\overline{11},\overline{10},10,10,\overline{9},8,8,\overline{7},7,\overline{5},\overline{3},\overline{2})$, we have $LN(\pi)=10$, $SN(\pi)=7$, $\ell_{O\leq N}(\pi)=6$ and $\ell_{O\geq N}(\pi)=5$. Moreover,  there are two parts $\overline{7},7$ of size $7$ in $\pi$, then we replace $\overline{7}$ in $\pi$ by ${7}$ to get $\lambda=(\overline{13},\overline{11},\overline{10},10,10,\overline{9},{8},8,7,7,\overline{5},\overline{3},\overline{2})$. Clearly, we have  $\ell_{O\leq N}(\lambda)=\ell_{O\leq N}(\pi)-1=5$.

\end{itemize}

For the four examples above, if we apply the map $\phi$ to the resulting overpartition $\lambda$, then we can recover the overpartition $\pi$.

\end{document}